\documentclass[letterpaper]{article}

\usepackage{anysize}
\usepackage[margin=1in]{geometry}
\usepackage[latin1]{inputenc}
\usepackage{amsfonts}
\usepackage{amssymb}
\usepackage{amsthm}
\usepackage{amsmath}
\usepackage{newlfont}

\begin{document}

\everymath{\displaystyle}

\begin{center}
{\huge Sum of Cubes is Square of Sum}\\
\vspace{3mm}
{ Edward Barbeau \footnote{barbeau@math.toronto.edu}, Samer Seraj \footnote{samer.seraj@mail.utoronto.ca}}\\
{ University of Toronto}\\
\vspace{3mm}
\end{center}

\textbf{Abstract.} Inspired by the fact that the sum of the cubes of the first $n$ naturals is equal to the square of their sum, we explore, for each $n$, the Diophantine equation representing all non-trivial sets of $n$ integers with this property. We find definite answers to the standard question of infinitude of the solutions as well as several other surprising results.

\section{Introduction}

An intriguing mathematical fact is that for every natural $n$,
\[1^3+2^3+\cdots+n^3=(1+2+\cdots+n)^2.\]

Thus the set $\{1,2,\ldots,n\}$ has the property that the sum of the cubes of its elements is equal to the square of the sum of its elements. It is natural to enquire about other sets $\{a_1,a_2,\ldots, a_n\}$ where this phenomenon occurs \[a_1^3+a_2^3+\cdots+a_n^3=(a_1+a_2+\cdots+a_n)^2.\]

The purpose of this paper is to present a number of new results and to encourage the reader to continue the investigation of this beautiful Diophantine equation.\\

We will call such a set a CS-set, and denote it by angular brackets $\langle a_1,a_2,\ldots, a_n\rangle$. Repeated elements are allowed but to avoid trivialities we will exclude sets that contain $0$ or contain both $k$ and $-k$ for some integer $k$. When all entries are positive, we refer to a positive CS-set; a CS$n$-set will denote a CS-set with $n$ elements.

\section{Positive CS-sets}

The following is an extraordinary generalization attributed to the French mathematician Liouville [2], [5].\\

\textbf{Proposition 1.} For each natural $n$, denote by $\tau(n)$ the number of positive divisors of $n$. Then $\langle \tau(d):d|n, d\ge 1\rangle$ is a CS-set, where $\tau$ is taken over the positive divisors of $n$.\\

\textit{Proof.} The proposition is clearly true for prime powers since for each prime $p$ and each natural $n$, the corresponding set of $p^{n-1}$ is
\[\langle\tau(p^0), \tau(p^1), \ldots, \tau(p^{n-1})\rangle=\langle 1,2,\ldots,n\rangle.\]
$\tau$ is well-known to be a multiplicative function [7], meaning that if $a,b$ are relatively prime integers then $\tau(ab)=\tau(a)\tau(b)$. By elementary multiplicative number theory, the functions
\[f(n)=\sum_{d|n}{\tau(d)} \hspace{3mm}\text{and}\hspace{3mm}F(n)=\sum_{d|n}{[\tau(d)}]^3\]
must also be multiplicative. All we need now is that $F(n)=[f(n)]^2$. Since this is already true for prime powers, the rest follows from the prime factorization of $n$ and the multiplicative property of $f$ and $F$. \hfill $\square$\\

\textbf{Proposition 2.} For each natural $n$, there are finitely many positive CS$n$-sets.\\

\textit{Proof.} Suppose that $a_k$ are the elements of a CS$n$-set such that $m$ is the largest entry. Then
\[m^3 \le \sum_{k=1}^n a_k^3 = \left( \sum_{k=1}^n a_k \right)^2 \le (nm)^2 = n^2m^2,\]
so that $m\le n^2$. Thus the entries in a positive CS$n$-set cannot exceed $n^2$, leaving finitely many $n$-tuples. \hfill $\square$\\

There is a striking general property of positive CS$n$-sets.\\

\textbf{Proposition 3.} For every natural $n$, there is precisely one positive CS$n$-set with distinct elements.\\

\textit{Proof.} We prove by induction that:\\

\textit{For integers $a_k$, if $1 \leq a_1 < a_2 < \cdots < a_n$ then
\[a_1^3 + a_2^3 + \cdots + a_n^3 \ge (a_1 + a_2 + \cdots + a_n)^2,\]
with equality if and only if $a_k = k$ for $1 \le k \leq n$.}\\

This is clear for $n = 1$. Assume that the proposition holds for some natural $n$; we will prove it for $n+1$. Note that $a_k \ge k$ and that $a_{n+1} - k \ge a_{n+1-k}$ for all values of $k$. We begin with the fact that
\[(a_{n+1} - n)(a_{n+1} - n - 1) \ge 0\]
with equality if and only if $a_{n+1} = n+1$. Then expanding gives
\begin{align*}
	a_{n+1}^2 - a_{n+1} &\ge 2na_{n+1} - n(n+1) = 2 \left[ na_{n+1} - {{n(n+1)}\over 2} \right]\\
	&= 2 \sum_{k=1}^n (a_{n+1} - k) \ge 2 \sum_{k=1}^n a_{n+1-k} = 2 \sum_{k=1}^n a_k
\end{align*}
with equality if and only if $a_k = k$ for all $k$. By the induction hypothesis, we have that
\begin{align*}
	\sum_{k=1}^{n+1} a_k^3 &= a_{n+1}^3 + \sum_{k=1}^n a_k^3 \ge a_{n+1}^2 + 2a_{n+1}\sum_{k=1}^n a_k + \left(\sum_{k=1}^n a_k \right)^2\\
	&= \left(a_{n+1} + \sum_{k=1}^n a_k \right)^2 = \left( \sum_{k=1}^{n+1} a_k \right)^2
\end{align*}
as desired. \hfill $\square$\\

If we do not require all the members of the set to be distinct, we see, by iterating through all possible $n$-tuples with entries within the upper bound from Proposition 2, that CS-sets occur surprisingly often.\\

$n = 2: \langle 1, 2 \rangle$, $\langle 2, 2 \rangle$\\
$n = 3: \langle 1, 2, 3 \rangle$, $\langle 3, 3, 3 \rangle$\\
$n = 4: \langle 1, 2, 2, 4 \rangle$, $\langle 1, 2, 3, 4 \rangle$, $\langle 2, 2, 4, 4\rangle$, $\langle 4, 4, 4, 4 \rangle$\\
$n = 5: \langle 1, 2, 2, 3, 5 \rangle$, $\langle 1, 2, 3, 4, 5 \rangle$, $\langle 3, 3, 3, 3, 6 \rangle$, $\langle 3, 3, 3, 4, 6 \rangle$, $\langle 5, 5, 5, 5, 5 \rangle$.\\

At this point, the number of possibilities increases markedly. For $n$ equal to $6$, $7$ and $8$ entries there are, respectively, $18$, $30$ and $94$ positive CS$n$-sets.\\

Mason [6] has indicated why CS$n$-sets are so frequent, by showing how they can be derived from arbitrary sets of $n$ natural numbers. For two sets $\langle a_i \rangle$ and $\langle b_j \rangle$ of integers, he defines their \textit{bag product} as the set $\langle a_i b_j \rangle$ of all products with one element from each set. The bag product of two CS-sets is also a CS-set. Mason constructs positive CS-sets from the bag product of an arbitrary set and a suitable constant set.\\

By taking bag products of sets of the form $\langle 1, 2, \ldots , n \rangle$, we can obtain Liouville's Tau generalization. More precisely, it can be seen from the proof of Proposition 1 that the Liouville set of $m$ is the bag product of the CS-sets corresponding to the largest power of each prime dividing $m$.\\

When $m = 24=2^3\cdot 3$, for example, this yields the set
\[ \langle \tau (1), \tau (2), \tau (3), \tau (4), \tau (6), \tau (8), \tau (12), \tau (24) \rangle = \langle 1, 2, 2, 3, 4, 4, 6, 8 \rangle,\]
which is the bag product of $\langle 1, 2, 3, 4 \rangle$ (corresponding to $2^3$) and $\langle 1, 2 \rangle$ (corresponding to $3$).\\

A full characterization of all positive CS-sets or all CS-sets in general may be out of reach, but we can improve the search efficiency. For example, we showed that the upper bound of an element of a positive CS$n$-set is $n^2$, but numerical evidence indicates that a much better bound can be established. Good bounds on the number of positive CS$n$-sets for each $n$ seems plausible as well.

\section{Extending CS-sets}

Numerical evidence shows that there are many examples of CS-sets that can be extended by a single entry to give another CS-set. We examine when this occurs.\\

\textbf{Proposition 4.} A CS-set can be extended by an element if and only if its sum $a$ is equal to $z(z-1)/2$ for some integer $z$. Moreover the appended integer is $z$. If $-z$ is a member of the CS-set, then we can replace the appending action with deleting $-z$.\\

\textit{Proof.} Suppose that we have a CS-set for which the sum of the entries is $a$; then the sum of the cubes of the entries is $a^2$. Suppose that appending an integer $z$ to the set results in another CS-set such that the new sum is $b$ and the new cube sum is $b^2$. Then
\[ \begin{cases} a+z=b \\ a^2+z^3=b^2 \end{cases} \implies \begin{cases} b-a=z \\ b+a=z^2 \end{cases} \implies \begin{cases}a=  z(z-1)/2 \\ b = z(z+1)/2 \end{cases}.\]

Note that the appended integer is $b-a=z$. Conversely, it is straightforward to show that if the sum of a CS-set is $z(z-1)/2$ for some integer $z$ then appending $z$ yields another CS-set. \hfill $\square$\\

For example, extending $\langle 1 \rangle$ yields in turn all of the sets of the form $\langle 1, 2, \ldots , n \rangle$. However, there are other possibilities, the smallest of the positive CS-sets being $\langle 1, 1, 4, 5, 5, 5 \rangle$ to which the numbers $\{ 7, 8, 9, \ldots\}$ can be appended in turn. There are many others but it is unknown for each $n$ how many CS$n$-sets exist with sum equal to $z(z-1)/2$ for some integer $z$.\\

Observe that for $a\ne 0$ the quadratic equation $a=z(z-1)/2\implies z^2-z-2a=0$ in $z$ has both a positive and a negative root. Thus we can alter the set by appending a negative integer $z$ instead. For example, beginning with $\langle 6, 6, 6, 6, 6, 6 \rangle$ gives
\[ \langle -8, 6, 6, 6, 6, 6, 6 \rangle \longrightarrow \langle -8, -7, 6, 6, 6, 6, 6, 6 \rangle \longrightarrow \langle -8, -7, 6, 6, 6, 6, 6 \rangle.\]

Overall a way of finding CS-sets is to start with any CS-set whose sum is a triangular number. Then we can follow the appending process of obtain a chain of CS-sets.\\

We can carry out a similar procedure to characterize when two numbers $x$ and $y$ can be appended to a suitable CS-set to produce another CS-set. If we take $y$ to be the negative of one of the numbers in the set, then we can simply replace the $-y$ in the set with $x$.\\

\textbf{Proposition 5.} A CS-set with sum $a$ can be extended by two entries to produce another CS-set if and only if $2(2a+1)=u^2+v^2+(u+v)^2$ for some integers $u,v$.\\

\textit{Proof.} Suppose that we have a CS-set with sum $a$ and cube sum $a^2$, and that appending some pair of integers $x$ and $y$ makes the sum $b$ and cube sum $b^2$. Then
\[ \begin{cases} b - a = x + y \\ b^2 - a^2 = x^3 + y^3 \end{cases} \implies b + a = x^2 - xy + y^2. \]
It follows that
\[\begin{cases}
	2a = (x^2 - xy + y^2) - (x + y) \implies 2(2a+1) = (y - x)^2 + (x - 1)^2 + (y - 1)^2\\
	2b = (x^2 - xy + y^2) + (x + y) \implies 2(2b+1) = (y - x)^2 + (x + 1)^2 + (y + 1)^2
\end{cases}.\]

Conversely, suppose the sum of a CS-set is $a$ and $2(2a+1)=u^2+v^2+(u+v)^2$ for some integers $u,v$. Let
\begin{align*}
	\begin{cases} x=u+1\\ y=u+v+1 \end{cases} &\implies \begin{cases} u=x-1 \\ v=y-x \end{cases} \implies 2a=(x^2-xy+y^2)-(x+y)\\
	&\implies 2(a+x+y)=(x^2-xy+y^2)+(x+y).
\end{align*}
Then we have that
\begin{eqnarray*}
	4(a^2+x^3+y^3)&=&[(x^2-xy+y^2)-(x+y)]^2+4(x^3+y^3)\\
	&=& [(x^2-xy+y^2)+(x+y)]^2\\
	&=& [2(a+x+y)]^2\\
	&=& 4(a+x+y)^2.
\end{eqnarray*}
Thus $a^2+x^3+y^3=(a+x+y)^2$. The process can be continued indefinitely since
\[2(2(a+x+y)+1)=(y-x)^2+(x+1)^2+(y+1)^2.\]
\hfill $\square$\\

For example, consider the set $\langle 3, 3, 3, 3, 4, 6, 8 \rangle$ whose sum $a$ is $30$ and whose cube sum is $30^2$. Since $2(2a + 1) = 122 = 4^2 + (-9)^2 + (-5)^2$, we can let $(x, y) = (5, -4)$ and so obtain the set $\langle 3, 3, 3, 3, 5, 6, 8 \rangle$.

\section{Zero Sums}

There is significant motivation to explore CS-sets with zero sum. We note that a constant integer multiple of the elements of a zero-sum CS-set is a zero-sum CS-set, as is the union of zero-sum CS-sets, and the union of a CS-set with any CS-set. Recall that the bag product of two CS-sets is also a CS-set, but the bag product of a zero-sum CS-set with any set is a CS-set.\\

\textbf{Proposition 6.} For $n=1,2,3,4$, there are no zero-sum CS$n$-sets.\\

\textit{Proof.} If $n=1$, we get the CS-set consisting of just $0$, which is excluded.\\

The $n=2$ case gives
\[a_1^3+a_2^3=0=(a_1+a_2)^2\implies a_1=-a_2,\]
which is excluded.\\

Note that for $n\ge 3$, if $\langle a_1,a_2,\ldots, a_n\rangle$ is a CS-set such that
\[a_1^3+a_2^3+\cdots+a_n^3=0=(a_1+a_2+\cdots+a_n)^2\]
then we have that
\[a_1^3+a_2^3+\cdots+a_{n-1}^3=(a_1+a_2+\cdots+a_{n-1})^3,\]
which is an interesting Diophantine equation in its own right and deserves attention. Its resolution would fully solve the zero-sum problem.\\

For $n=3$, this is
\[a_1^3+a_2^3=(a_1+a_2)^3\implies a_1 a_2(a_1+a_2)=0,\]
The equation holds when $a_1=0$ or $a_2=0$ or $a_1=-a_2$, all of which are excluded.\\

For $n=4$, this is
\[a_1^3+a_2^3+a_3^3=(a_1+a_2+a_3)^3\implies 0=(a_1+a_2+a_3)^3-(a_1^3+a_2^3+a_3^3)=3(a_1+a_2)(a_2+a_3)(a_3+a_1).\]
The solutions are $a_1=-a_2$ or $a_2=-a_3$ or $a_3=-a_1$, all of which are excluded. \hfill $\square$\\

A special family for $n=5$ can be found as follows.\\

\textbf{Proposition 7.} There are infinitely many CS$5$-sets with zero sum such that the entries share no positive factor greater than $1$.\\

\textit{Proof.} Inspired by the zero-sum CS-set $\langle -8,-7,1,5,9\rangle$, we try to obtain zero-sum CS$5$-sets of the form \[\langle -x,-y,r-s,r,r+s\rangle.\]
We require that
\[x+y=3r\hspace{3mm}\text{and}\hspace{3mm}x^3+y^3=(r-s)^3+r^3+(r+s)^3=3r(r^2+2s^2).\]
Subsequently,
\begin{eqnarray*}
	(3r)^3 &=& (x+y)^3\\
	&=& x^3+y^3+3xy(x+y)\\
	&=& 3r(r^2+2s^2) +3xy(3r)\\
	\implies 9r^2 &=& r^2+2s^2+3xy\\
	\implies (3x)\cdot(3y) &=& 24r^2-6s^2.
\end{eqnarray*}

We know that $3x+3y=9r$, so $3x$ and $3y$ are the roots of the quadratic equation
\[z^2-9rz +(24r^2-6s ^2)=0\implies\frac{z}{3}=\frac{9r\pm\sqrt{3(8s^2-5r^2)}}{6}.\]
Suppose that $8s^2 - 5r^2 = 3$. Then $r$ is odd and so $\frac{z}{3} = \frac{3r\pm 1}{2}$ yield integers $x, y$. Letting $t = 4s$, we get $t^2 - 10r^2 = 6$.\\

A solution of this Pellian equation is $(t, s) = (4, 1)$. Since the fundamental solution of $\alpha^2 - 10 \beta^2 = 1$ is $(\alpha, \beta) = (19, 6)$, a set of solutions of $t^2 - 10r^2 = 6$ is given by
\[t_k + r_k \sqrt{10} = (4 + \sqrt{10})(19 + 6 \sqrt{10})^k\]
where $k$ is any positive integer. For all of these solutions, $t$ is a multiple of $4$ so that $s$ is indeed an integer.\\

Generally, we get the CS$5$-sets
\[\left\langle-\frac{3r+1}{2}, -\frac{3r-1}{2}, r-s, r, r+s\right\rangle\]
where $8s^2 - 5r^2 = 3$. Note that the first two entries are apart by $1$ so the only positive common divisor of all the entries is $1$.\\

However, we can get other possibilities by taking the discriminant of the quadratic to be the square of any multiple of $3$. There are surely many CS$5$-sets that do not include an arithmetic progression. \hfill $\square$\\

It turns out that larger CS-sets with zero sum occur abundantly. This fact can be related to the Tarry-Escott problem [2]. This old problem asks for pairs of integer sets with the same cardinality and equal $k^{\text{th}}$ power sums for all $k$ running from $1$ to some natural $m$. We denote this property for two such sets $A, B$ by $A\overset{m}=B $.\\

\textbf{Proposition 8.} For $n=7,8$ and each $n\ge 10$, there are infinitely many zero-sum CS$n$-sets.\\

\textit{Proof.} Let $A=\{a_1, a_2,\ldots, a_n\}$ and $B= \{b_1, b_2,\ldots, b_n\}$. Frolov [4] states that for every integer $c$,
\begin{eqnarray*}
	\{a_1, a_2,\ldots, a_n\} &\overset{m}=& \{b_1, b_2,\ldots, b_n\}\\
	\implies \{c+a_1, c+a_2,\ldots, c+a_n\} &\overset{m}=& \{c+b_1, c+b_2,\ldots, c+b_n\}.
\end{eqnarray*}
This is easily verified upon expansion by the binomial theorem. We will concentrate on $m=3$. Suppose for integers $a_k$ and $b_k$ that $\sum_{k=1}^{n}{a_k^2}=\sum_{k=1}^{n}{b_k^2}$. We will call a pair of sets of equal cardinality and equal sum of squares a square pair, or an SP. Like CS$n$-sets, we will say SP$n$ to refer to a square pair of sets of cardinality $n$. It can be seen that
\[\{-a_k : 1 \le k \le n \}\cup\{a_k : 1\le k\le n\}\overset{3}=\{ -b_k : 1\le k \le n \}\cup\{b_k : 1 \le k \le n \},\]
since the sum and cube sum are both $0$. Then we have for any integer $c$ that
\[\{c-a_k : 1 \le k \le n \}\cup\{c+a_k : 1\le k\le n\}\overset{3}=\{c-b_k : 1\le k \le n \}\cup\{c+b_k : 1 \le k \le n \},\]
so both the sum and cube sum of each set
\[\langle\{c-a_k : 1 \le k \le n \}\cup\{c+a_k : 1\le k\le n\}\cup\{-c+b_k : 1\le k \le n \}\cup\{-c-b_k : 1 \le k \le n \}\rangle\]
is $0$, as desired.\\

Before expressing the full power of this technique, let us digress with an example. Begin with any Pythagorean triple, say $\{3, 4, 5\}$. Since $3^2+4^2=0^2+5^2$,
\begin{eqnarray*}
	\{ -3, -4, 3, 4 \}&\overset{3}=&\{ 0,-5,0,5 \}\\
	\implies \{ c-3, c-4, c+3, c+4 \}&\overset{3}=&\{c, c-5, c, c+5 \}.
\end{eqnarray*}
This gives rise to the family of CS-sets
\[\langle c-3, c-4, c+3, c+4,-c,-c+5,-c,-c-5 \rangle.\]
Taking $c = 6$ yields
\[\langle 3,2,9,10,-6,-1,-6,-11 \rangle = \langle -11,-6,-6,-1,2,3,9,10 \rangle.\]

By setting $c$ to the negative of an element, an entry can be strategically eliminated. For example, setting $c=3$ above yields
\[\langle 0,-1,6,7,-3,2,-3,-8\rangle=\langle -8,-3,-3,-1,2,6,7\rangle.\]
It is well known that there are infinitely many Pythagorean triples [7], so we can use them in this manner to get infinitely many CS$n$ sets with zero-sum for $n=7,8$.\\

We can get infinitely many SP$n$ for each $n$ by extending the idea of Pythagorean triples to Pythagorean $n$-tuples. For $n\ge 4$, define a Pythagorean $n$-tuple as an $n$-element list $(x_1,x_2,\ldots,x_n)$ for which $\sum_{k=1}^{n-1}{x_k^2}=x_n^2$. Since $\{x_1,\ldots,x_{n-1}\}$ and $\{x_n,0,\ldots,0\}$ form an SP, we want infinitely many Pythagorean $n$-tuples for each $n\ge 4$.\\

It holds for all integers $a,b$ that \[(a^2-b)^2+(2a)^2 b=(a^2+b)^2.\] Suppose that $n\ge 4$ and let $a$ and $\{a_1,a_2,\ldots,a_{n-2}\}$ be arbitrarily chosen integers. Let $b=\sum_{k=1}^{n-2}{a_k^2}$. Then \[(a^2-b)^2+\sum_{k=1}^{n-2}{(2aa_k)^2}=(a^2+b)^2,\]
which means that \[(a^2-b,2aa_1,\ldots, 2aa_{n-2},a^2+b)\] is a Pythagorean $n$-tuple. We can arrange that all of the $a_k$ are distinct or that several of them are equal.\\

Infinitely many Pythagorean $(n+1)$-tuples leads to infinitely many CS$(4n)$-sets for $n\ge 3$. An entry can be strategically eliminated by setting the translation $c$ to be the negative of a distinct element, so infinitude also holds for $4n-1$ elements when $n\ge 3$.\\

Since some of the initial $a_k$ can be chosen to be equal, two or three entries can be arranged to be equal. Those two or three entries can be annihilated by choosing $c$ to be their negative. So there are also infinitely many zero-sum CS$(4n-2)$-sets for $n\ge 3$, and CS$(4n-3)$-sets for $n\ge 4$. Thus we have shown that for $n=7,8$ and each $n\ge 10$ there are infinitely many CS$n$-sets with zero-sum. \hfill $\square$\\

Here is how a Pythagorean $4$-tuple can be used. Let $r$, $s$ and $c$ be arbitrary integers. Then \[(s^2-2r^2,2rs,2rs,s^2+2r^2)\]
is a Pythagorean $4$-tuple. As above, we can use it to construct an infinite family of CS-sets:
\begin{align*}
	\langle\{c-s^2+2r^2,c-2rs,c-2rs\}&\cup\{c+s^2-2r^2,c+2rs,c+2rs\}\\
	&\cup\{-c,-c,-c+2r^2+s^2\}\cup\{-c,-c,-c-2r^2-s^2\}\rangle.
\end{align*}
By our choice of $c$, we can obtain a CS-set of length $10$, $11$ or $12$.\\

If we want more variety among the entries, we can combine two disjoint Pythagorean $n$-tuples $A=(a_1, a_2,\ldots, a_n)$ and $B= (b_1, b_2,\ldots, b_n)$ to get an SP$n$ with distinct elements:
\[\sum_{k=1}^{n-1}{a_k^2}=a_n^2, \hspace{3mm}\sum_{k=1}^{n-1}{b_k^2}=b_n^2\implies b_n^2+\sum_{k=1}^{n-1}{a_k^2}=a_n^2+\sum_{k=1}^{n-1}{b_k^2}.\]

The cases that remain are $n=6,9$. There are solutions such as \[\langle-11,-5,-4,2,8,10\rangle\hspace{3mm} \text{and}\hspace{3mm} \langle-17,-10,-8,-1,2,3,6,7,18\rangle,\] but apart from taking constant multiples, we have not succeeded in describing non-trivial infinite families.

\section{Infinitude}

\textbf{Proposition 9.} For $n=1,2$, there are finitely many CS$n$-sets.\\

\textit{Proof.} The only CS$1$-set is $\langle 1 \rangle$, since $a^3=a^2\implies a=1$.\\

In the $n=2$ case, for $a\ne b$, we have
\begin{eqnarray*}
	&& a^3 + b^3 = (a + b)^2\\
	\implies && a^2 -ab + b^2 = a + b\\
	\implies && a^2 -(b + 1)a + (b^2 - b) = 0.
\end{eqnarray*}
The discriminant is
\[(b+1)^2-4(b^2-b)=-3b^2+6b+1,\] which has finitely many positive values. Iterating through them, the only CS$2$-sets are $\langle 1, 2 \rangle$ and $\langle 2, 2 \rangle$. \hfill $\square$\\

Suppose that we begin with any set of $n$ integers $\{ a_1, a_2, \ldots , a_n \}$. Let
\begin{eqnarray*}
	&&u = a_1^3 + a_2^3 + \cdots + a_n^3,\\
	&&v = a_1 + a_2 + \cdots + a_n.
\end{eqnarray*}

Can we get a CS-set from this by multiplying each entry by an integer $t$? It is straightforward to check that $t$ should be $v^2/u$. Thus, we can construct a CS-set whenever $v^2$ is a multiple of $u$. This will always occur when $u = 1$, $u = 2$ and when $u = v$.\\

For example, any solution of the equation $x^3 + y^3 + z^3 = 1$ will generate a CS$3$-set, and it was shown by Ramanujan [2] that there are infinitely many such non-trivial triples.\\

We proceed with the $u=v$ case for $n = 3$ using a method attributed to Chowla.\\

\textbf{Proposition 10.} For $n=3,4$, there are infinitely many CS$n$-sets with distinct elements.\\

\textit{Proof.} From numerical examples, we notice that there seems to be many cases where
\[x^3+y^3+z^3=x+y+z, \hspace{3mm} x+y=3q, \hspace{3mm} z=-2q\]
for some integer $q$. Suppose that this is true for some integers $x,y,z$. Then
\begin{eqnarray*}
	q&=& x+y+z\\
	&=& x^3+y^3+z^3\\
	&=& (x+y)[(x+y)^2-3xy]+z^3\\
	&=& 3q[(3q)^2-3xy]+(-2q)^3\\
	&=& 19q^3-9qxy\\
	\implies xy&=& \frac{19q^2-1}{9}.
\end{eqnarray*}

As a result,
\begin{eqnarray*}
	x+y&=&3q\\
	\implies (x-y)^2&=&9q^2-4xy\\
	&=& 9q^2-4\left(\frac{19q^2-1}{9}\right)\\
	&=&\frac{5q^2+4}{9}\\
	\implies [3(x-y)]^2-5q^2&=&4.
\end{eqnarray*}

Thus there is motivation to explore the Pellian equation $\alpha^2 - 5 \beta^2 = 4$, which has infinitely many solutions $(\alpha_k, \beta_k)$ given by $\alpha_k + \beta_k \sqrt{5} = (3 + \sqrt{5})(9 + 4\sqrt{5})^k$ where $k$ is a non-negative integer. Then
\begin{eqnarray*}
	\alpha_{k+1} &=& 9 \alpha_k + 20 \beta_k\\
	\beta_{k+1} &=& 4\alpha_k + 9 \beta_k.
\end{eqnarray*}
Since $\alpha_{k+1} \equiv -\beta_k$ and $\beta_{k+1} \equiv \alpha_k \pmod 3 $, we see that $\alpha_k$ is a multiple of 3 whenever $k$ is even. Therefore, there are infinitely many pairs $(\phi, \psi)$ of positive integers for which $9\phi^2 - 5\psi^2 = 4$.\\

Noting that $\phi$ and $\psi$ are necessarily of the same parity, for each such solution let
\[(x,y,z) = \left( {{\phi + 3\psi}\over 2}, {{3\psi - \phi}\over 2}, -2\psi \right).\]
It is straightforward to check that the sum and cube sum of $x, y, z$ are both $\psi$.\\

For $n = 4$, simply append to this triple an entry equal to $1$ or $-1$. In this construction, the CS$n$-sets for $n=3,4$ have distinct elements.\hfill $\square$\\

Note that \[a^3+b^3+c^3=a+b+c\iff \binom{a+1}{3}+\binom{b+1}{3}=\binom{-c+1}{3},\] assuming $a,b>0>c$. Thus the problem is asking for pairs of tetrahedral numbers which sum to another tetrahedral number
\[\binom{x}{3}+\binom{y}{3}=\binom{z}{3}.\]

This problem was attacked by several mathematicians. Segal [9] proved for the first time in $1962$ that the only solution with $x=y$ is $x=y=4$ and $z=5$. He also noted computational efforts before his time. Chowla's proof of infinitude is repeatedly mentioned in the literature [3] [8] [9], though it seems that his original paper was not published or has been lost; it does not appear in his collected works. Since he is mentioned in Segal's paper, it can be assumed that he found the proof before $1962$. His manipulations were extended by Edgar [3] in $1964$ to find two further families of solutions. It is curious that Edgar refers to a paper of Chowla with Newman, Segal and Wunderlich which is ``to appear''. In $1965$, Oppenheim [8] published a new method of finding infinitely many solutions.\\

At long last, we complete our discussion with a twist that the reader may have foreseen. We apply a past proposition to easily deduce the strongest result.\\

\textbf{Proposition 11.} For each $n\ge 5$ that there are infinitely many CS$n$-sets with distinct elements.\\

\textit{Proof.} For a fixed $n\ge 5$, we simply append to $\langle 1,2,\ldots, n-5\rangle$ a zero-sum CS$5$-set with all entries distinct and modulus greater than $n-5$; there are certainly infinitely many such CS$5$-sets by the constructions in Proposition 7 and their constant multiples. \hfill $\square$

\section{Conclusion}

In our journey, we have answered the most pressing questions about CS-sets. The collected results are:

\begin{enumerate}

	\item For each natural $n$:
	\begin{enumerate}
		\item[a.] There are finitely many positive CS$n$-sets.
		\item[b.] There is precisely one positive CS$n$-set with distinct entries, namely $\langle 1,2,\ldots, n\rangle$.
	\end{enumerate}
	
	\item For any CS-set with sum $a$:
	\begin{enumerate}
		\item[a.] It can be extended by an entry to produce another CS-set if and only if for some integer $z$, \[a=z(z-1)/2.\]
		\item[b.] It can be extended by two entries to produce another CS-set if and only if for some integers $u,v$, \[2(2a+1)=u^2+v^2+(u+v)^2.\]
	\end{enumerate}
	
	\item
	\begin{enumerate}
		\item[a.] For $n=1,2,3,4$, there are no CS$n$-sets with sum zero.
		\item[b.] For each natural $n\ge 5$, there are infinitely many CS$n$-sets with sum zero.
	\end{enumerate}
	
	\item
	\begin{enumerate}
		\item[a.] For $n=1,2$, there are respectively $1,2$ total CS$n$-sets.
		\item[b.] For each natural $n\ge 3$, there are infinitely many CS$n$-sets with distinct entries.
	\end{enumerate}

\end{enumerate}

The few results that we have given on CS-sets probably only scratch the surface. Mentioned throughout are remaining related open problems of interest. There are undoubtedly many other connections to be made.\\

\section{References}

\begin{enumerate}

\item Edward J. Barbeau, \textit{Pell's Equation}. Springer, 2003.

\item Edward J. Barbeau, \textit{Power Play}. Mathematical Association of America, 1997. pp. 5, 14-15, 33

\item Hugh Maxwell Edgar, Some remarks on the Diophantine equation $x^3 + y^3 + z^3 = x + y + z$, \textit{Proceedings of the American Mathematical Society} \textbf{16} (February 1965) 148-153

\item M. Frolov, Égalités à deux degrés, \textit{Bulletin de la Société Mathématique de France} \textbf{17} (1889) 69-83

\item Ross Honsberger, \textit{Ingenuity in Mathematics}, Mathematical Association of America, 1970

\item John Mason, Generalizing ``Sums of cubes equals to squares of sums'', \textit{The Mathematical Gazette} \textbf{85} (March 2001) 50-58

\item Ivan Niven, Herbert S. Zuckerman, Hugh L. Montgomery, \textit{An Introduction to the Theory of Numbers}, John Wiley and Sons, Inc., 1991. pp. 190, 231

\item Alexander Oppenheim, On the Diophantine equation $x^3 + y^3 + z^3 = x + y + z$, \textit{Proceedings of the American Mathematical Society} 17 (1966), 493-496

\item S. L. Segal, A Note on Pyramidal Numbers, \textit{The American Mathematical Monthly} \textbf{69} (1962) 637-638

\end{enumerate}

\end{document}